\documentclass[11pt]{article}
\usepackage{psfig}
\usepackage{latexsym}
\usepackage{graphicx}
\usepackage{psfrag}
\newcommand {\be}[1]{\begin{equation}\label{#1}}
\newcommand {\ee}{\end{equation}}
\newcommand {\bea}{\begin{eqnarray}}
\newcommand {\eea}{\end{eqnarray}}
\newcommand {\vvert}{\par\vspace{.1in}}

\newcommand{\proof}{{\bf Proof. }}
\newcommand{\qed}{\hfill $\Box$}

\newtheorem{theorem}{Theorem}[section]
\newtheorem{lemma}[theorem]{Lemma}

\newtheorem{coro}[theorem]{Corollary}

\newcommand{\pr}{{\mbox{\rm Prob}}}

\newcommand{\x}{{\bf x}}
\newcommand{\y}{{\bf y}}
\newcommand{\z}{{\bf z}}

\newcommand{\vv}{{\bf v}}
\newcommand{\zprime}{{\bf z'}}
\newcommand{\vprime}{{\bf v'}}
\newcommand{\X}{{\bf X}}

\newcommand{\N}{{\bf N}}
\newcommand{\0}{{\bf 0}}
\begin{document}
\title{The diameter of a  long range percolation graph}
\author{ Don Coppersmith
\thanks{IBM T.J. Watson Research Center, Yorktown Heights, NY 10598, USA.
        Email address: dcopper@us.ibm.com }
\and
David Gamarnik
\thanks{IBM T.J. Watson Research Center, Yorktown Heights, NY 10598, USA.
        Email address: gamarnik@watson.ibm.com }
\and
Maxim Sviridenko
\thanks{IBM T.J. Watson Research Center, Yorktown Heights, NY 10598, USA.
        Email address: sviri@us.ibm.com }
}
\date{}
\maketitle
\begin{abstract}
We consider the following long range percolation model:
an undirected  graph with the node set $\{0,1,\ldots,N\}^d$,
has edges $(\x,\y)$  selected with probability $\approx \beta/||\x-\y||^s$
if $||\x-\y||>1$,
and with probability $1$ if $||\x-\y||=1$, for some parameters $\beta,s>0$.
This model was introduced by Benjamini and Berger \cite{benjamini_berger}, who
obtained bounds on the diameter of this graph for the one-dimensional case $d=1$ 
and for various values
of $s$, but left cases $s=1,2$ open.
We show that, with high probability,
the diameter of this graph is  $\Theta(\log N/\log\log N)$ when $s=d$,
and, for some constants $0<\eta_1<\eta_2<1$, 
it is at most $N^{\eta_2}$, when $s=2d$ and is
at least $N^{\eta_1}$
when $d=1,s=2,\beta<1$ or $s>2d$.
We also provide a simple proof that the diameter is at most $\log^{O(1)}N$
with high probability, when  $d<s<2d$,
established previously in \cite{benjamini_berger}.
\end{abstract}

\section{Introduction}
Long-range percolation is a model in which
any two elements $x,y$ of some (finite or countable) metric space
are connected by edges with some probability, inverse proportional
to the distance  between the points. The motivation for studying
this model is dual. First, it   naturally extends  a classical
percolation models on a lattice, by adding edges between non-adjacent nodes
with some positive probability.
The questions of existence of infinite components were considered specifically 
by Schulman \cite{schulman},
Aizenman and Newman \cite{aizenman_newman} and Newman and Schulman
\cite{newman_schulman}, where the metric space is $\mathcal{Z}$
and edges $(i,j)\in \mathcal{Z}^2$ are selected with probability
$\beta/|i-j|^s$ for some parameters $\beta,s$. Existence of such an infinite component
with positive probability usually implies its existence
with probability one, by appealing to Kolmogorov's $0-1$ law.
It was shown in \cite{newman_schulman} and in \cite{aizenman_newman} respectively,
that percolation occurs if $s=2,\beta>1$ and (suitably defined) short range
probability is high enough,  and does not occur if $s=2,\beta\leq 1$, for
any value of the short range probability.
 
The second motivation for studying long-range percolation is
modeling social networks, initiated by Watts and Strogatz \cite{watts_strogatz}.
They considered a random graph model on  integer points of a circle, in which neighboring nodes
are always connected by an edge, and, in addition, each node is connected
to a constant number of other nodes uniformly chosen from a circle.
Their motivation was a famous experiment conducted by Milgram \cite{milgram}, which essentially studied
the diameter of the ``social acquaintances'' network and
introduced the notion of ``six degrees of separation''. Watts and Strogatz
argued that their graph provides a good model for different types of networks,
not only social networks (world wide web, power grids), and showed that the diameter
of their random graph is much smaller than the size of the graph. This model
was elaborated later by Kleinberg \cite{kleinberg_smallw},
who considered a model similar to a long-range percolation
model on a two-dimensional grid, although the work was concerned mostly
with algorithmic questions of constructing simple
decentralized algorithms for finding short paths between the nodes.

The present paper is motivated by a recent work
by Benjamini and Berger \cite{benjamini_berger}. They consider a
one-dimensional long-range percolation model
in which the nodes are elements of a finite circle $\{0,1,\ldots,N\}$.
An edge $(i,j)$ exists with probability
one if $|i-j|=1$, and with probability $1-\exp(-\beta/|i-j|^s)$ otherwise,
for some parameters $\beta,s$, here the distance $|\cdot|$ is taken
with respect to a circle. Since for large $|i-j|$,
$1-\exp(-\beta/|i-j|^s)\approx\beta/|i-j|^s$, this model is closely
related to the infinite percolation model on $\mathcal{Z}$, with
an important distinction, however. The graph is finite and, since neighboring nodes
are connected with probability one, the graph is connected. Thus,  the
percolation question is irrelevant as such; rather, as in models of ``social networks'',
the diameter of the graph is of interest.
It is shown in \cite{benjamini_berger} that the diameter of the circle
graph above is with high probability is a
constant, when $s<1$;
is $O(\log^{\delta}N)$,
for some $\delta>1$, when $1<s<2$; and is linear  ($\Theta(N)$), when $s>2$.
These results apply immediately to a graph on an interval $\{0,1,\ldots,N\}$. 
A multidimensional version  of this problem
with a graph on a node set $\{0,1,\ldots,N\}^d$ was also considered
by Benjamini, et.al in \cite{bkps}, who showed that the diameter is $\lceil d/(d-s)\rceil$
when $s<d$. 
The cases $s=1,2$ were left open in \cite{benjamini_berger}
and the authors conjectured that the diameter is $\Theta(\log N)$
when $s=1$, and $\Theta(N^{\eta})$ for some $0<\eta<1$, when $s=2$.
In addition, the authors conjectured that, for  the
case $1<s<2$, $\log^{\delta} N$ is also a lower bound for some $\delta>1$. In other words,
the system experiences a phase transition at $s=1$ and $s=2$.

In this work we consider a multidimensional version of the problem.
Our graph  has a node set $\{0,1,\ldots,N\}^d$ and edges are selected 
randomly using a long-range percolation  $\beta/||\x-\y||^s$ law. 
We obtain
upper and lower bounds on the diameter for the regimes $s=d,d<s<2d,s=2d$ and
$s>2d$. This corresponds to regimes $s=1,1<s<2,s=2,s>2$ for the one-dimensional case. 
We show that,
with high probability, for $s=d$, the diameter of this graph is
$\Theta(\log N/\log\log N)$; for $d<s<2d$ the diameter is at most
$\log^{\delta} N$ for some constant $\delta>1$;
and for $s=2d$, the diameter
is at most $N^{\eta_2}$, for some constant $0<\eta_2<1$. We also prove a
lower bound $N^{\eta_1},\eta_1<1$ on the diameter, which holds with high probability
but only when $d\geq 1,s>2d$ or  $d=1,s=2,\beta<1$. We do not have lower bounds 
for other cases. Note that our  lower bound for $s>2d$ is weaker than known linear
lower bound when $d=1$. We conjecture that the linear lower bound holds for general 
dimensions.
Our results, when applied to the one-dimensional case, 
support bounds conjectured  in \cite{benjamini_berger}
for the case $s=2$ and disprove it for the case $s=1$. 
As we mentioned above, the upper bound $\log^{\delta}N$ for the case $d<s<2d$
was proven in \cite{benjamini_berger} for the one-dimensional case.
It was pointed to the authors, that the  proof extends to a multidimensional 
case as well.
We  provide here an alternative proof which seems simpler. 
Summarizing the results of present paper and of \cite{benjamini_berger}, the diameter
of the long-range percolation graph in one-dimensional case 
experiences a phase transition at $s=1,2$ and
has a qualitatively different values for $s<1;\, s=1;\, 1<s<2;\, s=2$ and $\beta<1;\, s>2$.
Whether the same holds true for for general dimensions (whether $s=d,s=2d$ are the
only critical values) remains to be seen. Our results only partially support this
conjecture.

\section{Model and the main result}
Our model is a random graph $G=G(N)$ on a node set $[N]_d\equiv\{0,1,\ldots,N\}^d$ -
integral points of the $d$-dimensional cube with side length $N$.
Let $||\x||$ denote an ${\bf L_1}$ norm in the space ${\cal Z}^d$.
That is $||\x||=\sum_{i=1}^dx_i$.
Nodes $\x,\y\in [N]_d$ are connected with probability $1$ if $||\x-\y||=1$,
and, otherwise, with probability $1-\exp(-{\beta\over ||\x-\y||^s})$, where $\beta>0, s>0$ are
some fixed parameters. Let $D(N)$ denote the (random) diameter of the graph $G(N)$,
and let $P(N)$ denote the (random) length of a shortest path between nodes 
$\0\equiv (0,\ldots,0)$ and $\N=(N,\ldots,N)$.
For any $\x,\y\in [N]_d$ let also $P(\x,\y)$ denote the length of a shortest path
between  nodes $\x,\y$ in the graph $G(N)$.
Our main result is as follows.
 
\begin{theorem}\label{main_result}
There exist constants $C_1,C_2,C_s>0,\delta>1,\,\,0<\eta_1<\eta_2<1$,
which in general depend on $s$, $\beta$ and on dimension $d$  
such that
\begin{enumerate}
\item $\lim_{N\rightarrow\infty}
\pr\{D(N)\geq N^{\psi}\}=1$, for any $s>2d,\psi<{s-2d\over s-d-1}$.
\item $\lim_{N\rightarrow\infty}\pr\{D(N)\leq N^{\eta_2}\}=1$, for $s=2d$
and \\
$\lim_{N\rightarrow\infty}\pr\{D(N)\geq N^{\eta_1}\}=1$, for $d=1,s=2,\beta<1$.
\item $\lim_{N\rightarrow\infty}\pr\{C_s\log N\leq D(N)\leq \log^{\delta}(N)\}=1$, for $d<s<2d$.
\item $\lim_{N\rightarrow\infty}
\pr\{{C_1\log(N)\over\log\log N}\leq D(N)\leq {C_2\log(N)\over\log\log N}\}=1$, for $s=d$.
\end{enumerate}
\end{theorem}

As we mentioned above,  it was shown in \cite{bkps} that the diameter is, with high probability,
$\lceil d/(d-s)\rceil$, when $s<d$. Also 
part 3 of the theorem above was proven
by Benjamini and Berger in \cite{benjamini_berger} for the one-dimensional case. 
They also pointed out to the authors that their proof holds for a multidimensional case as well.
We provide here
a simpler proof.
Throughout the paper we use standard notations $f=O(g),f=\Omega(g),f=\Theta(g),f=o(g)$,
which mean respectively that for two functions $f(N),g(N)$,
$f(N)\leq C_1g(N),f(N)\geq C_2g(N),C_3g(N)\leq f(N)\leq C_4g(N),
f(N)/g(N)\rightarrow 0$, for some constants $C_i,i=1,2,3,4$ which in general
depend on $\beta,s$, but do not depend on $N$.
Also, throughout the paper
$[n]_d$ denotes an integral cube $\{0,1,\ldots,n\}^d$ for
any nonnegative  integer $n$. The logarithmic function is always
assumed to be with the base $e$.

\section{Case $s>2d$. Lower bound}
In this section we show that,with high probability, the diameter of the graph $G(N)$ is
at least  essentially $N^{s-2d\over s-d-1}$.
As we noted, for the one dimensional case $d=1$ this is weaker than
the existing linear lower bound $\Omega(N)$ (\cite{benjamini_berger}).
\vvert
 
{\it Proof of Theorem \ref{main_result}, Part 1:}
We fix a constant $\psi<{s-2d\over s-d-1}$. For any
$k>N^{1-\psi}$ let $L(k)$ be the total
number of edges between pairs of points at distance
exactly $k$. We will now show that if $\psi<(s-2d)/(s-d-1)$ then
$\sum_{k>N^{1-\psi}}kL(k)\leq dN/2$, with high probability. 
Since $||\N||=dN$, then
this would imply that, with high probability, any path between $\0$ and $\N$ 
would contain at least $dN/(2N^{1-\psi})=(dN^{\psi})/2$ edges 
and the proof would be completed.
For a fixed pair of nodes $\x,\y$ at a distance $k$, the probability that
the edge between them exist is $1-\exp(-\beta/k^s)\leq \beta/k^s$,
where we use $\exp(-\beta x)\geq 1-\beta x$ for all $0\leq x\leq 1$.
For a fixed node $\x$ there are $\Theta(k^{d-1})$ nodes $\y$
which are at distance $k$ from $\x$; also there are $N^d$ choices 
for the node $\x$. Combining $E[L(k)]= O(N^dk^{d-1}(\beta/k^s))$.
Then
\[
\sum_{k>N^{1-\psi}}kE[L(k)]= 
O(\beta N^d\sum_{k>N^{1-\psi}}k^{d-s})=O(N^dN^{(1-\psi)(d-s+1)}).
\]
For the given choice of $\psi$, we have $d+(1-\psi)(d-s+1)<1$
and the value above is $o(N)$.
Using Markov's inequality, we obtain 
\[
\pr\{\sum_{k>N^{1-\psi}}kL(k)>N/2\}\leq {o(N)\over (N/2)}=o(1). 
\]
\qed

\section{Case $s=2d$.}
 
\subsection{Upper bound}
In this subsection we prove that when $s=2d$, there exists
a constant $0<\eta<1$, which depends on $\beta$ and $d$, 
such that with high probability
$D(N)\leq N^{\eta}$. To this end we
first establish an upper bound on $\max_{\x,\y\in [N]_d}E[P(\x,\y)]$ and then
use this bound to obtain a polynomially small bound on $\pr\{D(N)> N^{\eta}\}$
for some constant $\eta<1$.

\vvert
 
{\it Proof of Theorem \ref{main_result}, Part 2:}

We first assume that $N$ is a power of $3: N=3^m$, for some
integer $m>0$, and then consider the general case.
For any fixed integer $n$ let $R(n)=\max_{\x,\y\in [n]_d}E[P(\x,\y)]$.
That is, $R(n)$ is the maximum over expected lengths of shortest paths between all the
pairs of points in the cube $[n]_d$. We obtain an upper bound
on $R(N)$ by relating it to $R(N/3)$. Divide the cube $[N]_d$
into $3^d$ subcubes of the type 
$I_{i_1\ldots i_d}\equiv \prod_{j=1}^d[i_j{N\over 3},(i_j+1){N\over 3}],
0\leq i_j\leq 2$. Each cube has a side length $N/3$ (which is integer since 
$N$ is a power of three). We say that two
such cubes are neighboring if they have at least a common  node.
For example $[0,N/3]^d$ and $[N/3,2N/3]^d$ are neighboring through a corner node
$(N/3,\ldots,N/3)$. We now fix  a pair of points $\x,\y\in [N]_d$ and estimate
$P(\x,\y)$ by considering two cases.

\begin{enumerate}
\item $\x,\y$ belong to the same subcube $I=I_{i_1\ldots i_d}$. The length of a
shortest path between these two points using edges of $[N]_d$
is not bigger than the length of the shortest path between same points but
using only edges of the subcube $I$. Therefore $E[P(\x,\y)]\leq R(N/3)$.

\item $\x,\y$ belong to different   subcubes $I,I'$. 
Let $E=E(I,I')$ be the event ``there exists at least one
edge between some nodes $\vv\in I,\vprime\in I'$. The probability
that $E$ occurs is at least $1-\exp(-\beta ({N\over 3}+1)^{2d}/(dN)^{2d})$
since there are $({N\over 3}+1)^d$ nodes in each cube,  and the largest
possible distance between them is $dN$. In particular, $\pr\{E\}$ is not
smaller than a certain constant $\delta>0$, independent of $N$. 
We now estimate $E[P(\x,\y)]$ conditioned on $E$ and $\bar E$.
Given that $E$ occurs, select an edge $(\vv,\vprime)$ between the cubes $I,I'$.
Then 
\[
E[P(\x,\y)|E]\leq E[P(\x,\vv)|E]+E[P(\vprime,\y)|E]+1
\]
Note, however, that edges within each cube $I,I'$ are selected 
independently from edges between cubes and specifically are independent
from the event $E$. Therefore, since $\x,\vv$ belong the the same cube,
$E[P(\x,\vv)|E]\leq R(N/3)$. Similarly, $E[P(\vprime,\y)|E]\leq R(N/3)$.
We conclude $E[P(\x,\y)|E]\leq 2R(N/3)+1$.
Now, suppose $E$ does not occur. Select a cube $I''$ which is a 
neighboring cube for cubes $I,I'$ (it is easy to see that such
a cube exists). Specifically, let $\z (\zprime)$ be the nodes shared by cubes
$I$ and $I''$ ($I'$ and $I''$). Then arguing as above
$E[P(\x,\y)|\bar E]\leq E[P(\x,\z)|\bar E]+E[P(\z,\zprime)|\bar E]+E[P(\zprime,\y)|\bar E]
\leq 3R(N/3)$. Combining, we obtain
\[
E[P(\x,\y)]\leq (2R(N/3)+1)\pr\{E\}+3R(N/3)(1-\pr\{E\})=(3-\pr\{E\})R(N/3)+\pr\{E\}\leq
\]
\[
(3-\delta)R(N/3)+1.
\]
\end{enumerate}
We conclude, $R(N)=\max_{\x,\y\in [N]_d}E[P(\x,\y)]\leq (3-\delta)R(N/3)+1$.
Applying this bound $m-1=\log N/\log 3-1$ times, we obtain
\[
R(N)\leq (3-\delta)^{m-1}R(3)+\sum_{i=0}^{m-2}(3-\delta)^i=O((3-\delta)^m)=
O(N^{\log(3-\delta)\over \log 3}),
\]
Note, $\alpha\equiv\log(3-\delta)/\log 3<1$. We obtain
$R(N)=O(N^{\alpha})$ for some $\alpha<1$.

In order to generalize the bound for all $N$, it is tempting to argue
that $R(N)\leq R(3^{m})$ as long as $N\leq 3^{m}$.
This would require proving a seemingly obvious
statement that $R(n)$ is a non-decreasing function of $n$.
While this is most likely correct, proving it does not
seem to be trivial. Instead, we proceed as follows.
Let $m$ be such that $3^m\leq N<3^{m+1}$. We cover the
cube $[N]_d$ with $3^d$ cubes $I_i, i=1,\ldots,3^d$ with side length $3^m$,
with a possible overlapping. Specifically, $I_i\subset [N]_d$ and
$\cup_i I_i=[N]_d$. Let $\x,\y\in [N]_d$ be arbitrary. Find
cubes $I_{i_1},I_{i_2},I_{i_3}$ such that $\x\in I_{i_1},\y\in I_{i_3}$
and $I_{i_1}\cap I_{i_2}\neq\emptyset, I_{i_2}\cap I_{i_3}\neq\emptyset$.
Let $\z_1,\z_2$ be some nodes lying in these intersections. Then
$E[P(\x,\y)]\leq E[P(\x,\z_1)]+E[P(\z_1,\z_2)]+E[P(\z_2,\y)]=O((3^m)^{\alpha})$,
where the last equality follows since 
pairs $(\x,\z_1),(\z_1,\z_2),(\z_2,\y)$ lie within  cubes $I_{i_1},I_{i_2},I_{i_3}$
respectively and each of them has a side length $3^m$. But $3^m\leq N$. We conclude
$E[P(\x,\y)]=O(N^{\alpha})$ and $R(N)=\max_{\x,\y}E[P(\x,\y)]=O(N^{\alpha})$.

We now finish the proof of part 2, upper bound, by obtaining 
a similar bound on the  diameter $D(N)$.
Fix an arbitrary $0<\epsilon,\gamma<1$ such that $\alpha+\epsilon<1$ and
$\epsilon-d(1-\gamma)>0$.
Divide the cube $[N]_d$ into  equal subcubes
$I_{i_1\ldots i_d}=\prod_{j=1}^d [i_jN^{\gamma},(i_j+1)N^{\gamma}],\,\,
0\leq i_j\leq N^{1-\gamma}$ 
each with side length $N^{\gamma}$. The total number of subcubes is $N^{d(1-\gamma)}$.
Fix any such cube $I$ and let
$\x(I)$ be its lower corner (the node with smallest possible coordinates).
We showed above $E[P(\0,\x(I))]\leq O(N^{\alpha})$, from which, using Markov inequality, 
\[
\pr\{P(\0,\x(I))\geq N^{\alpha+\epsilon}\}=O({N^{\alpha}\over N^{\alpha+\epsilon}})=
O({1\over N^{\epsilon}}).
\]
Then
\[
\pr\{\max_I P(\0,\x(I))\geq N^{\alpha+\epsilon}\}=
O({N^{d(1-\gamma)}\over N^{\epsilon}})=O({1\over N^{\epsilon-d(1-\gamma)}}).
\]
On the other hand for every cube $I$ and every $\x\in I$ we have
trivially, $P(\x,\x(I))\leq dN^{\gamma}$. 
Since $D(N)\leq 2\sup_{\x\in [N]_d}P(0,\x)$, then
\[
\pr\{D(N)\geq 2(dN^{\gamma}+N^{\alpha+\epsilon})\}=
O({1\over N^{\epsilon-d(1-\gamma)}})=o(1).
\]
We take $\eta=\max\{\gamma,\alpha+\epsilon\}<1$ and obtain
$\pr\{D(N)\geq 4dN^{\eta}\}=o(1)$.
This completes the proof of the upper bound. \qed

\subsection{Lower bound}
The proof of the lower bound for the one-dimensional case $d=1,s=2,\beta<1$ is similar to the
proof for the case $s>2$, from  \cite{benjamini_berger}
and uses the notion of a cut point.
We first show that $E[D(N)]\geq N^{\eta}$ for a certain
constant $0<\eta<1$, for large $N$. Then we show that this bound holds
actually with high probability.
Given a node $1\leq i\leq N-1$, we call it a cut node if
there are no edges which go across $i$. Namely, $i$
is a cut point if edges $(j,k)$ do not exist for all $j<i<k$.
The probability that $i$ is a cut node is
$\exp(-\beta\sum_{j<i<k}{1\over |j-k|^2})\geq
\exp(-\beta\sum_{1\leq n\leq N}{n-1\over n^2})=\Theta({1\over N^{\beta}})$.
Then the expected number of cuts is $\Omega(N^{1-\beta})$
(which will be helpful  to us only if $\beta<1$).
But the shortest path $P(N)$ and as a result the diameter $D(N)$ are
not smaller  than the number of cuts.
Taking $\eta<1-\beta$, we obtain the bound $E[D(N)]\geq N^{\eta}$
for large $N$.
 
To prove that the bound actually holds with high probability, we
need the following lemma. A node $i\in [N]$ is defined to be an isolated
node if it is only connected to $i-1,i+1$ (or only to $1$ or $N-1$ if $i=0$
or $N$ respectively). We denote by $I(N)$ the number of isolated nodes in $[N]$.

\begin{lemma}\label{isolation}
For any $\beta>0$,
the number of isolated nodes $I(N)$ is $\Omega(N)$ with probability tending
to one as $N\rightarrow\infty$.
\end{lemma}
 
\proof Given $i\in [N]$, the probability that it is isolated is
\[
p(i)\equiv\exp(-\beta\sum_{j\leq N} 1/|i-j|^2)\geq\exp(-2\beta\sum_{j\leq N}1/j^2)\geq
\exp(-\beta\pi^2/3)\equiv c(\beta).
\]
where we use the fact $\sum_{j=1}^{\infty}1/j^2=\pi^2/6$.
Then the expected number of isolated points is $E[I(N)]=\sum_ip(i)\geq c(\beta)N$.
We now estimate
$E[I^2(N)]=\sum_ip(i)+\sum_{i \neq j}\pr\{i,j\,\,{\rm are\,\, isolated}\}$.
But
\[
\pr\{i,j\,\,{\rm are\,\, isolated}\}=\pr\{j\,\,{\rm isolated}|i\,\,{\rm isolated}\}p(i)=
{p(j)p(i)\over \exp(-{\beta\over |i-j|^2})}=p(j)p(i)\exp({\beta\over |i-j|^2}).
\]
Then
\[
E[I^2(N)]=\sum_ip(i)+\sum_{i \neq j}p(j)p(i)\exp({\beta\over |i-j|^2})=
\sum_ip(i)+\sum_{|i-j|<\sqrt{N}}p(j)p(i)\exp({\beta\over |i-j|^2})+
\]
\[
\sum_{|i-j|\geq\sqrt{N}}p(j)p(i)\exp({\beta\over |i-j|^2}).
\]
Note that $\exp({\beta\over |i-j|^2})=1+O(1/N)$ when
$|i-j|\geq \sqrt{N}$. Also $\exp({\beta\over |i-j|^2})\leq \exp(\beta)=O(1)$
for all $i,j$, and $\sum_{|i-j|\geq\sqrt{N}}p(j)p(i)\leq
\sum_{i,j}p(j)p(i)=E^2[I(N)]$.
Combining,
\[
E[I^2(N)]\leq \sum_ip(i)+\sum_{|i-j|<\sqrt{N}}O(p(j)p(i))+
\sum_{i,j}(1+O({1 \over N}))p(j)p(i) \leq O(N)+O(N^{3\over 2})+
(1+O({1 \over N}))E^2[I(N)].
\]
Since $N^2 \geq E^2[I(N)]\geq c^2(\beta)N^2$ then
${\rm Var}(I(N))=O(N^{3/2}) = O(1/\sqrt{N})E^2[I(N)]$.
Using Chebyshev's inequality,
\[
\pr\{I(N)<(1/2)c(\beta)N\}\leq \pr\{I(N)\leq (1/2)E[I(N)]\}\leq
{{\rm Var}(I(N))\over (1/4)E^2[I(N)]}=O(1/\sqrt{N}).
\]
We proved that the number of isolated points is linear in $N$ with high probability. \qed
 
\vvert
 
We now complete the proof of the lower bound. Divide the interval
$[N]$ into $N^{2\over 3}$ intervals $I_1,I_2,\ldots,I_{N^{2\over 3}}$
each of length $N^{1\over 3}$. Construct a new graph on $N^{2\over 3}$
nodes, such that for $1\leq i,j\leq N^{2\over 3}$, and $|i-j|>1$, nodes $i,j$ are connected
by an edge if and only if there exists at least one edge between intervals $I_i,I_j$
in the original graph, and, for $|i-j|=1$, the nodes $i,j$ are
always connected. If $|i-j|>1$, then the probability
that nodes $i,j$ of the new graph are not connected is at least
$\exp(-\beta N^{2\over 3}/((|i-j|-1)^2N^{2\over 3})\geq \exp(-\beta'/|i-j|^2)$
for $\beta'=4\beta$. Applying Lemma \ref{isolation}, the new graph
has $\Omega(N^{2\over 3})$ isolated nodes with high probability. In other words, the
original graph has $\Omega(N^{2\over 3})$ intervals $I_i$ which do not have
edges connecting it to other intervals except for possibly intervals
$I_{i-1},I_{i+1}$. We call these intervals isolated intervals
and the number of them is denoted by $N_{iso}=\Omega(N^{2\over 3})$.
For each interval $I_i$ and each $x\in I_i$, we
say that $x$ is a local cut point if it is a cut point with
respect to just the graph induced by vertices from  $I_i$. We showed above that the expected
number of local cut points is at least $|I_i|^{\eta}=N^{\eta\over 3}$,
for any $\eta<1-\beta$ and for all $i$. Let $C(I_i)$ be the number of local cut points
in the interval $I_i$.  We now show that, with high probability,
at least one of the isolated intervals has at least $(1/2)N^{\eta\over 3}$ 
local cut points.
Note $C(I_i)$ are independent from each other.
We have $E[C(I_i)]\geq N^{\eta\over 3}$. Also ${\rm Var}(C(I_i))\leq |I_i|^2=N^{2\over 3}$.
Applying Chebyshev's inequality, we have
\[
\pr\{\sum_{I_i{\rm\,\, isolated}}C(I_i)/N_{iso}<
{1\over 2}N^{\eta\over 3}\}\leq 
{{\rm Var(I_i)}\over {1\over 2}N^{\eta\over 3}N_{iso}}=O({1\over N^{\eta\over 3}}),
\]
where the last equality follows from $N_{iso}=\Omega(N^{\eta\over 3})$.
Therefore, with high probability, at least one of the isolated intervals contains 
at least $(1/2)N^{\eta\over 3}$ local cut points.
We denote this interval by $I_{i^*}$. Let us estimate the number of edges
between $I_{i^*}$ and $I_{i^*-1}$ and $I_{i^*+1}$. Note that in defining
isolated interval $I_{i^*}$ with many local cut points,
we only considered edges between intervals $I_i,I_j$
which are not adjacent to each other (to obtain an isolated interval)
and edges within intervals $I_i$ (when considering local cuts). Specifically,
we did not consider edges between adjacent intervals $I_{i^*},I_{i^*-1},I_{i^*+1}$.
Then, by independence, the expected number
of edges between $I_{i^*}$ and $I_{i^*-1}$ is at most
\[
\sum_{k=1}^{N^{1\over 3}}(k+1)(1-\exp(-{\beta\over k^2}))+O(1)=
O(\sum_{k=1}^{N^{1\over 3}}{\beta\over k})=O(\log N),
\]
where we use $\exp(-\beta x)\geq 1-\beta x$ for all $x\in [0,1]$.
Using Markov's inequality, the probability that the number of
edges between $I_{i^*}$ and $I_{i^*-1}$ is bigger than $\log^2N$
is at most $O(1/\log N)$. A similar bound applies to $I_{i^*+1}$.
We conclude that with high probability there are at most
$2\log^2N$ edges outgoing from interval $I_{i^*}$ totally
(as it is isolated, it can only be connected to intervals $I_{i^*-1},I_{i^*+1}$).
Since the number of local cuts in $I_{i^*}$ is $\Omega(N^{\eta\over 3})$
then there are two local cuts $i_1,i_2$, such interval $[i_1,i_2]$
contains at least $\Omega(N^{\eta\over 3}/\log^2N)=\Omega(N^{\eta\over 4})$
local cuts and no outside edges are connected to nodes in interval
$[i_1,i_2]$. Let the number of local cuts in $[i_1,i_2]$ be $L$.
We take the $(1/3)L$-th and the $(2/3)L$-th local cut in this interval.
By construction, the shortest path between these local cuts is at
least $(1/3)L=\Omega(N^{\eta\over 4})$. We conclude, $D(N)=\Omega(N^{\eta\over 4})$,
with high probability. \qed

\section{Case $d<s<2d$.}
The lower bound $D(N)\geq C_s\log N$ was proven 
to hold with high probability
in \cite{benjamini_berger} for the case $d=1$,  using
branching theory and the fact that for each node, the expected
number of its neighbors is a constant. The proof extends easily to all dimensions $d$.
We now focus on an upper bound. 
Our proof is similar to the one in \cite{benjamini_berger}
and is based on renormalization technique, although our analysis
is simpler.
 
\vvert
{\it Proof of Theorem \ref{main_result}, Part 3:}
We have $d<s<2d$. Let us fix $\alpha<1$ such that $2d\alpha>s$.
Split the cube $[N]_d$ into  equal subcubes 
$I_{i_1\ldots i_d}\equiv 
\prod_{j=1}^d[i_j\lceil N^{\alpha}\rceil,(i_j+1)\lceil N^{\alpha}\rceil-1]$
with side length $\lceil N^{\alpha}\rceil$. If $N/\lceil N^{\alpha}\rceil$
is not an integer then we make the cubes containing nodes $(\ldots,N,\ldots)$
overlap partially with some other cubes. In the following we drop
the rounding $\lceil\cdot\rceil$ for simplicity, the argument still holds.
Consider the following
event $E_1$: ``there exist two cubes $I,I'$ such that
no edge exists between points $\x\in I$ and $\y\in I'$''.
Each resulting cube $I=I_{i_1\ldots i_d}$ we split further into subcubes
with side length $N^{\alpha^2}$. We
consider the event $E_2$: ``there exist a cube $I$ with side length $N^{\alpha}$ 
and its two subcubes $I_1,I_2$ with side length $N^{\alpha^2}$, such that no edge exists
between points in $I_1$ and $I_2$''. We continue this
process  $m$ times, obtaining in the end cubes with side
length $N^{\alpha^m}$. Assume that none of the events $E_1,E_2,\ldots,E_m$
occurs. We claim that then the diameter of our original graph is
at most $2^{m+1}N^{\alpha^m}$. In fact, since event $E_1$ does not
occur any two points $\x,\y\in [N]_d$ are connected by a path
with length at most $2\bar D(N^{\alpha})+1$, where
$\bar D(N^{\alpha})$ is the (random) largest diameter of
the cubes  $I_{i_1\ldots i_d}$ with side length $N^{\alpha}$. 
Similarly, since event $E_2$ does not occur,
$\bar D(N^{\alpha})\leq 2\bar D(N^{\alpha^2})+1$, where $\bar D(N^{\alpha^2})$
is the largest diameter of  the subcubes with side length $N^{\alpha^2}$, 
obtained in second stage. In the end we obtain
that the diameter of our graph satisfies $D(N)\leq 2^mD(N^{\alpha^m})+2^m\leq
2^{m+1}dN^{\alpha^m}$, since trivially, $D(N^{\alpha^m})\leq dN^{\alpha^m}$.
We now show that for a certain value of $m$ this upper bound on the
diameter $D(N)$ is at most $\log^{\delta} N$ for some constant $\delta>1$ and
simultaneously, the probability $\pr\{\wedge_{r=1}^m \bar E_r\}\rightarrow 1$,
as $N\rightarrow\infty$.
For a given cube with side length $N^{\alpha^{r-1}}$ and its two given
subcubes with side length $N^{\alpha^r}$, the probability that no
edges exist between these two subcubes is at most
$\exp(-\beta N^{2d\alpha^r}/(dN)^{s\alpha^{r-1}})=
\exp(-\Theta(N^{\alpha^{r-1}(2d\alpha-s)}))$, since
there are $N^{2d\alpha^r}$ pairs of points considered and the largest distance
among any two of them is $dN^{\alpha^{r-1}}$.
Since there are at most
$N^{2d}$ pairs of such subcubes, then the probability of the event $E_r$ is
bounded above by $N^{2d}\exp(-\Theta(N^{\alpha^{r-1}(2d\alpha-s)}))$.
We conclude
\[
\pr\{\vee_{r=1}^m E_r\}\leq \sum_{r=1}^mN^{2d}e^{-\Theta(N^{\alpha^{r-1}(2d\alpha-s)})}\leq
mN^{2d}e^{-\Theta(N^{\alpha^{m}(2d\alpha-s)})}
\]
Let us fix a large constant $C$ and take
\[
m={\log\log N-\log\log\log N+\log (2d\alpha-s)-\log C
\over\log {1\over\alpha}}=O(\log\log N).
\]
A straightforward computation shows that for this value of $m$,
\be{probbound}
\pr\{\vee_{r=1}^m\bar E_r\}=O(e^{-\Theta(\log^C N)}).
\ee
On the other hand, we showed above that, conditioned on event $\wedge_r\bar E_r$,
we have $D(N)=O(2^mdN^{\alpha^m})$. For our choice of $m$ a simple calculation
shows that $\alpha^m\log N=O(\log\log N)$ or $N^{\alpha^m}=\log^{O(1)}N$.
Also, since $m=O(\log\log N)$, then $2^m=O(\log^{O(1)}N)$. This completes the
proof. \qed
 
\vvert
 
In the course of the proof we established the following bound
which follows immediately from  (\ref{probbound}).
\begin{coro}\label{tail} For any constant $C$, there
exists a constant $\delta>1$ such that
\[
\pr\{D(N)>\log^{\delta}N\}\leq O(e^{-\Theta(\log^C N)})
\]
\end{coro}

\section{Case $s=d$.}\label{s1}

{\it Proof of Theorem \ref{main_result}, Part 4:}
We first prove a lower bound. We show that $D(N)\geq (d-\epsilon)\log N/\log\log N$
with high probability, for any constant $0<\epsilon<1$. Observe, that,
for any $1<k\leq N$ and for each node $\x\in [N]_d$, there are $\Theta(k^{d-1})$,
nodes at distance $k$ from $\x$. Each such node is connected to $\x$
with probability $1-\exp(-\beta/k^d)\leq \beta/k^d$.
(We used $\exp(-\beta x)\geq 1-\beta x$ for all $x\in [0,1]$).
Then the expected number of nodes connected to $\x$ by an edge is at most
$O(1)+O(\sum_{1\leq k\leq dN} (k^{d-1}/k^d))=O(\log N).$
Then, the total expected number of nodes which are reachable from $\x$ by
paths with length $\leq m$ is at most $c^m\log^mN$, for some constant $c$. 
We denote the number of such nodes $B(m)$. Using Markov's inequality
\[
\pr\{B(m)\geq N^d\}\leq {E[B(m)]\over N^d}\leq {c^m\log^mN\over N^d}\rightarrow 0
\]
if $m=(d-\epsilon)\log N/\log\log N$. Therefore, with probability tending
to one, the diameter $D(N)$ is $\Omega(\log N/\log\log N)$.
 
We now focus on a more difficult part -- the upper bound.
The proof is fairly technical, but is based on a simple observation
which we present now. We have already noted that
any fixed  node $\z$, in particular, node $\N$, has in expectation $\Theta(\log N)$ neighbors.
We will show later in the formal proof that this actually holds  with high probability.
Consider a subcube $I=[0,N/\log^c N]^d$ for a certain
constant $c$. Let $\y$ be a neighbor of $\x$. The probability that $\y$ has 
no neighbors  in $I$ is at most $\exp(-\beta N^d/(d^dN^d\log ^{cd}N))$, since the 
largest possible distance is $dN$ and the number of nodes in $I$ is
$N^d/\log^{cd}N$. Then probability that none of the $\Theta(\log N)$ neighbors
of $\N$ is connected to some node of $I$ by a path of length $\leq $ two
is at most $\exp(-\beta N^d\log N/(d^dN^d\log ^{cd}N))=\exp(-\beta\log^{1-cd}N)$.
If $c<1/d$ then this quantity converges to zero. Therefore, with high probability
$\N$ is connected to some node $\X_1\in I$ by a path of length $2$. Applying this 
argument for $\X_1$ we find a node $\X_2$ which is connected to $\X_1$ by
a path of length two and such that all the coordinates of $\X_2$ are at most
$N/\log^{2c} N$. Continuing
$m$ times we will obtain that $\N$ is 
connected by a path of length $O(m)$ to some node $\X_m$ with all the coordinates 
$\leq N/\log^{cm}N$. Taking $m=O(\log N/\log\log N)$ we will obtain
that, with high probability, $\N$ is connected to $\0$ by a path of length $\leq O(m)$.
We now formalize this intuitive argument.

We fix an arbitrary node $\z_0\in [N]_d$.
Consider all the  paths $(\x,\y,\z_0)$ with length two, which end in
node $\z_0$. That is  edges $(\x,\y),(\y,\z_0)$ exist. Let
$\X_1={\rm argmin} ||\x||$, where the minimum is taken over all such paths.
In other words, $\X_1$ is the smallest, in norm, node connected to $z_0$ via
a path of length at most $2$.
Note, $\X_1$ is random and $||\X_1||\leq ||\z_0||$,
as $\z_0$ is connected to itself by a path of length two.
Similarly, let $\X_2<\X_1$ be the smallest, in norm, node, connected to $\X_1$ via
a path of length   $2$. We continue this procedure for $m$ (to be defined later) steps
and obtain a (random) node $\X_m$.
 
\begin{lemma}\label{ym}
For any constantly large integer $c$, if
$m=(2d+2)\cdot 2^{c+1}\log N/\log\log N$, then
the bound  $||\X_m||\leq \exp((\log N)^{d/2^c})$ holds 
with probability at least $1-1/N^{2d}$.
\end{lemma}
 
Before we prove the lemma, let us show how it is used to
prove the result. We invoke part 3
of Theorem \ref{main_result}, which we proved in the previous
section.
Choose a constant integer $c$ such that
$2^c/d\geq 2\delta$, where $\delta>1$ is a constant from part 3
of Theorem \ref{main_result}. Applying  part 3 of Theorem \ref{main_result},
the diameter of the cube $[\exp((\log N)^{d/2^c})]_d$ is at most
$((\log N)^{d/2^c})^{\delta}\leq \log ^{1\over 2}N=o(\log N/\log\log N)$
with high probability. In particular
$\sup_{\x:||\x||\leq \exp((\log N)^{d/2^c})}P(0,\x)=o(\log N/\log\log N)$
with high probability.
By the conclusion of the lemma, with probability at least $1-O(1/N^{2d})$,
each fixed node $z_0\in [N]_d$  is connected to some node
$\X_m$ with $||\X_m||\leq\exp((\log N)^{d/2^c})$ by a path of length
$m=O(\log N/\log\log N)$. Then, with probability at least $1-O(1/N^d)$,
all the nodes $z_0\in [N]_d$ are connected to some corresponding nodes
$\X_m\in [\exp((\log N)^{d/2^c})]_d$ by a path of length
$O(\log N/\log\log N)$.
Combining, we obtain that $\sup_{\z_0\in [N]_d}P(0,\z_0)=O(\log N/\log\log N)$
with probability at least $1-o(1)$. But $D(N)\leq 2\sup_{\z_0\in [N]_d}P(0,\z_0)$. \qed
 
\vvert
 
{\it Proof of Lemma \ref{ym}:}
We fix a node $\x$ with  $||\x||\leq ||\z_0||$, fix  $1\leq r\leq m$ and 
consider $\X_r$ conditioned on event $\X_{r-1}=\x$ (assume $\X_0=\z_0$).
Our goal for the remaining part is the following 

\begin{lemma}\label{x}
If $||\x||>\exp((\log N)^{d\over 2^c})$, then
\be{rr-1}
E\Big[||\X_r||\Big | \X_{r-1}=\x\Big]\leq O({||\x||\over (\log N)^{1/2^{c+1}}}).
\ee
In other words, at each step $r=1,2,\ldots,m$, the expected value
of $||\X_r||$ decreases by a factor of $O({1\over (\log N)^{1/2^{c+1}}})$,
provided that $||\X_{r-1}||$ is still bigger than $\exp((\log N)^{d\over 2^c})$.
\end{lemma}

\proof Let $B(\x)$ be the total number of nodes 
 which are connected to $\X_{r-1}=\x$ and which have a norm smaller than $||\x||$. 
Note, that for each such node $\y$, $||\y-\x||\leq ||\y||+||\x||< 2||x||$.
We first show that with
probability at least $1-O({1\over (\log N)^{d/2^c}})$, the equality
$B(\x)=\Omega(\log ||\x||)$ holds. For any
fixed $k\leq 2||\x||$ there are $\Theta(k^{d-1})$ nodes $\y$
which for which $||\y-\x||=k$ and $||\y||<||\x||$. Each 
such node is connected by an edge to $\x$ with probability $1-\exp(-\beta/k^d)$.
Then 
\[
E[B(\x)]=\sum_{0\leq k<2||x||}(1-\exp(-{\beta\Theta(k^{d-1})\over k^d}))=\Theta(\log||x||)
\]
Let $c_1<c_2$ be constants, such that $c_1\log||x||\leq E[B(\x)]\leq c_2\log||x||$.
We now estimate the second moment
\[
E[B^2(\x)]=E[B(\x)]+
\sum_{\y_1\neq \y_2,||\y_1||,||\y_2||<||\x||}
(1-\exp(-{\beta\over ||\y_1-\x||}))(1-\exp(-{\beta\over ||\y_2-\x||}))
\leq
\]
\[
E[B(j)]+\sum_{||\y_1||,||\y_2||<||\x||}
(1-\exp(-{\beta\over ||\y_1-\x||}))(1-\exp(-{\beta\over ||\y_2-\x||}))=
\]
\[
E[B(\x)]+(E[B(x)])^2.
\]
It follows, ${\rm Var}(B(\x))\leq E[B(\x)]$.
Using Chebyshev's inequality, 
\[
\pr\{B(\x)\leq (1/2)c_1\log||\x||\}\leq
\pr\{|B(\x)-E[B(\x)]|\geq (1/2)c_1\log||\x||\}\leq
\]
\be{o1}
{{\rm Var}(B(\x))\over (1/4)c_1^2\log^2||\x||}\leq
{c_2\log||\x||\over (1/4)c_1^2\log^2||\x||}=O({1\over \log ||\x||})\leq
O({1\over(\log N)^{d/2^c}}),
\ee
where the last inequality follows from the assumption 
$||\x||>\exp((\log N)^{d/2^c})$ of the lemma.
Let 
\[
V(\x)=\{\z:||\z||\leq {||\x||\over(\log N)^{1/2^{c+1}}}\}
\]
In particular, $|V(\x)|=\Theta(||\x||^d/(\log N)^{d/2^{c+1}})$.
Suppose $\y,||\y||<||\x||$ is any node  connected by an edge to $\x$ (if any exist)
Note that the distance between $\y$ and any node in $V(\x)$
is smaller than $2||\x||$. Then, the probability
that $\y$ has no nodes  in $V(\x)$ connected to it
by an edge is at most
\[
\exp(-{\beta \Theta(||\x||^d)\over(\log N)^{d/2^{c+1}}||\x||^d})=
\exp(-{\Theta(1)\over(\log N)^{d/2^{c+1}}}).
\]
By (\ref{o1}), with  probability at least $1-O({d\over(\log N)^{1/2^c}})$,
$\x$ has   $\Omega(\log ||\x||)$
nodes $\y,||\y||<||\x||$  connected to it. Conditioned on this event,
the probability that no node in $V(\x)$
is connected to $\x$ by a path of length two is at most
$\exp(-{\Omega(\log ||\x||)\over (\log N)^{d/2^{c+1}}})$. By assumption,
$||\x||>\exp((\log N)^{d\over 2^c})$ or $\log ||\x||> (\log N)^{d\over 2^c}$,
using which, $\exp(-{\Omega(\log ||\x||)\over (\log N)^{d/2^{c+1}}})\leq
\exp(-{\Omega((\log N)^{d/2^{c+1}})})$.
It follows, that the probability that no node in $V(\x)$
is connected to $\x$ by a path of length two, is at most
\[
O({1\over(\log N)^{d/2^c}})+\exp(-{\Omega((\log N)^{d/2^{c+1}})})=
O({1\over(\log N)^{d/2^c}}).
\]
Summarizing, conditioned on $\X_{r-1}=\x$,
the bound $||\X_r||\leq {||\x||\over(\log N)^{1/2^{c+1}}}$ holds
with probability at least $1-O({d\over(\log N)^{1/2^c}})$.
On the other hand, with probability one $||\X_r||\leq ||\X_{r-1}||$. We conclude
\[
E\Big[||\X_r||\Big |\X_{r-1}=\x\Big]\leq {||\x||\over(\log N)^{1/2^{c+1}}}+
O({||\x||\over(\log N)^{d/2^c}})= O({||\x||\over(\log N)^{1/2^{c+1}}}).
\]
This completes the proof of Lemma \ref{x}. \qed
 
\vvert

We now complete the proof of  Lemma \ref{ym}. Note, that for any $2\leq r\leq m$,
$E[\X_r|\X_{r-1},\X_{r-2},\ldots,\X_1]=E[\X_r|\X_{r-1}]$. We denote $\exp((\log N)^{d/2^{c}})$
by $\alpha(N)$. We have,
\[
\pr\{||\X_m||>\alpha(N)\}=
\sum_{\alpha(N)<||\x_m||\leq ||\x_{m-1}||<||\z_0||}
\pr\{\X_m=\x_m|\X_{m-1}=\x_{m-1}\}\pr\{\X_{m-1}=\x_{m-1}\}\leq
\]
\[
\sum_{\alpha(N)<||\x_m||\leq ||\x_{m-1}||<||\z_0||}
||\x_m||\pr\{\X_m=\x_m|\X_{m-1}=\x_{m-1}\}\pr\{\X_{m-1}=\x_{m-1}\}\leq
\]
\[
\sum_{\alpha(N)<||\x_{m-1}||<||\z_0||}
E\Big[||\X_m||\Big|\X_{m-1}=\x_{m-1}\Big]\pr\{\X_{m-1}=\x_{m-1}\}.
\]
But, using bound (\ref{rr-1}) of Lemma \ref{x}, we have 
$E\Big[||\X_m||\Big|\X_{m-1}=\x_{m-1}]\leq O(||x_{m-1}||/(\log N)^{1/2^{c+1}})$,
as long as $||x_{m-1}||>\alpha(N)$. We obtain
\[
\pr\{||\X_m||>\alpha(N)\}\leq
O({1\over (\log N)^{1/2^{c+1}}})\sum_{\alpha(N)<||\x_{m-1}||<||\z_0||}
||\x_{m-1}||\pr\{\X_{m-1}=\x_{m-1}\}=
\]
\[
O({1\over (\log N)^{1/2^{c+1}}})\sum_{\alpha(N)<||\x_{m-1}||\leq ||\x_{m-2}||<||\z_0||}
||\x_{m-1}||\pr\{\X_{m-1}=\x_{m-1}|\X_{m-2}=\x_{m-2}\}\pr\{\X_{m-2}=\x_{m-2}\}\leq \]
\[
O({1\over (\log N)^{1/2^{c+1}}})
\sum_{\alpha(N)<||\x_{m-2}||<||\z_0||}
E\Big[||\X_{m-1}||\Big|\X_{m-2}=\x_{m-2}\Big]\pr\{\X_{m-2}=\x_{m-2}\}\leq
\]
\[
\Big(O({1\over (\log N)^{1/2^{c+1}}})\Big)^2
\sum_{\alpha(N)<||\x_{m-2}||<||\z_0||}||\x_{m-2}||\pr\{X_{m-2}=\x_{m-2}\},
\]
where in the last inequality we used bound (\ref{rr-1}) of Lemma \ref{x} again. 
Continuing this conditioning
argument $m-1$ times, we obtain that for some constant $C$
\[
\pr\{||\X_m||>\alpha(N)\}\leq {C^{m-1}\over (\log N)^{m-1\over 2^{c+1}}}||\z_0||\leq
{(\log N)^{1\over 2^{c+1}} C^m\over (\log N)^{m\over 2^{c+1}}}dN.
\]
But, by assumption of the lemma, $m=(2d+2)\cdot 2^{c+1}\log N/\log\log N$, from which
$(\log N)^{1\over 2^{c+1}}C^m=o(N)$ and $\pr\{||\X_m||>\alpha(N)\}\leq 1/N^{2d}$ 
for large $N$. \qed

\section{Concluding remarks and open questions}
We considered a long range percolation model on
an  graph with a node set $\{0,1,\ldots,N\}^d$.
Answering some open questions raised by Benjamini and Berger in
\cite{benjamini_berger}, we showed that if two nodes at a distance
$r$ are connected by an edge with probability $\approx\beta/r^s$,  then,
with high probability,
the diameter of this graph is $\Theta({\log N\over\log\log N})$ when $s=d$, and
is at most $N^{\eta}$  for some value $\eta<1$,
when $s=2d$. We also proved a lower bound $N^{\eta'},\eta'<1$ on the 
diameter for the cases $d=1,s=2,\beta<1$ and  $s>2d, d\geq 1$. 
Note that for the case $d=1,s>2$ our bound is weaker than known
linear lower bound $\Omega(N)$ established in \cite{benjamini_berger}.
We conjecture that this linear lower bound holds for all dimensions $d$
as long as $s>2d$. Other unanswered regimes are lower bounds for
$s=2d$; $d=1,s=2,\beta>1$; and $d<s<2d$. 
It would also be interesting to compute  the limits
${D(N)\over(\log N/\log\log N)}\rightarrow C$ and $\log D(N)/\log N\rightarrow \eta$
or even show that these limits actually exist when $s=d,2d$ respectively.

{\bf Acknowledgments.} We wish to thank I.Benjamini and N.Berger for clarifying
their work and identifying several errors in an earlier version of this paper.

\bibliographystyle{plain}
\bibliography{/.../watson.ibm.com/fs/users/G/gamarnik/My_Papers/bibliography.bib}

\end{document}